\documentclass[12pt]{amsart}
\usepackage{amscd,amsthm,amssymb,amsfonts,amsmath,epsfig,epic,bbm,url, stmaryrd, latexsym, tikz,mathrsfs, hyperref}
\usepackage[arrow,matrix,tips,2cell]{xy}
\usetikzlibrary{arrows}

\theoremstyle{plain}
\newtheorem{theorem}{Theorem}
\newtheorem{lemma}[theorem]{Lemma}

\theoremstyle{definition}
\newtheorem{definition}[theorem]{Definition}

\theoremstyle{remark}

\newtheorem*{acknowledgments}{Acknowledgments}

\newcommand{\Z}{\mathbb Z}    
\newcommand{\Q}{\mathbb Q}    
\newcommand{\R}{\mathbb R}    
\newcommand{\C}{\mathbb C}    
\newcommand{\T}{\mathbb T}    

\newcommand{\ZZ} {\mathbb{Z}}
\newcommand{\QQ} {\mathbb{Q}}
\newcommand{\RR} {\mathbb{R}}
\newcommand{\CC} {\mathbb{C}}

\newcommand{\TT}{\mathbb{T}}

\newcommand{\GL}{\operatorname{GL}}

\newcommand{\ignore}[1]{\relax}

\newcommand{\Hom}{\operatorname{Hom}}
\newcommand{\ra}{\rightarrow}

\newcommand{\bsd}{\operatorname{bsd}}

\newcommand {\Spec} {\operatorname{Spec}}

\newcommand {\tors}  {\operatorname{tors}}

\newcommand{\id}{\operatorname{id}}

\newcommand{\Conv}{\operatorname{Conv}}

\newcommand{\shS}{\mathcal{S}} 
\newcommand{\shX}{\mathcal{X}} 
\newcommand{\pp}{\mathcal{P}}

\renewcommand{\P}{\mathscr{P}}

\begin{document}

\title{Compactifying torus fibrations over integral affine manifolds with singularities}

\author{Helge Ruddat and Ilia Zharkov}
\address{Johannes Gutenberg-Universit\"at Mainz \& Universit\"at Hamburg}
\email{ruddat@uni-mainz.de,\ helge.ruddat@uni-hamburg.de}
\address{Kansas State University, 138 Cardwell Hall, Manhattan, KS 66506}
\email{zharkov@ksu.edu}

\begin{abstract}
This is an announcement of the following construction: given an integral affine manifold $B$ with singularities, we build a topological space $X$ which is a torus fibration over $B$.
The main new feature of the fibration $X\to B$ is that it has the discriminant in codimension~2. 
\end{abstract}
\maketitle

\section{Introduction}
\makeatletter
\def\blfootnote{\xdef\@thefnmark{}\@footnotetext}
\makeatother
\blfootnote{H.R. was supported by DFG grant RU 1629/4-1 and the Department of Mathematics at Universit\"at Hamburg.
The research of I.Z. was supported by Simons Collaboration grant A20-0125-001. }
There have been a lot of studies of half-dimensional torus fibrations and their integral affine structures on the base spaces inspired by the Strominger-Yau-Zaslow conjecture \cite{SYZ}. 
This area was very active in the beginning of the 2000's with many approaches of different flavor: topological \cite{Zh}, \cite{Gross}, symplectic \cite{Gross00}, \cite{Leung}, \cite{Rua}, \cite{J03},\cite{CBM}, \cite{Au07}, \cite{Au09}, \cite{Evans}, metric \cite{GW00}, \cite{KS00}, \cite{LYZ}, non-Archimedean \cite{KS04}, tropical \cite{Mikhalkin}, combinatorial \cite{HZ}, and log-geometric \cite{logmirror1},  \cite{logmirror2}, \cite{Pa07}. 
For surveys on the early developments, see \cite{T06,Gross09}. The toric case was considered in \cite{CL,CLL,FLTZ12}.
A more recent surge and interest is mostly tropical \cite{Mat,SS18,AGIS,Mi19,H19,MR}, non-Archimedean \cite{NXY} or topological \cite{P18}. For more recent surveys, we refer to \cite{Gross12,Ch}.
Broadly speaking, all this research developed into a new field of mathematics: tropical geometry.

In this note, we essentially follow the Gross-Siebert setup \cite{logmirror1},  \cite{logmirror2}, with some slight modifications. 
We replace the polyhedral decomposition of the base $B$ by a regular CW-decomposition for the gain of flexibility, cf.~the notion of ``symple'' in \cite{affinecoh}. Also we relax requirements for the monodromy by allowing arbitrary lattice simplices for local monodromies, not just the elementary ones. That requires a little more care for the local monodromy assumptions, but does not seem to affect the topological side of the story much. On the other hand, when we compare our model with the Kato-Nakayama space of a canonical Calabi-Yau family, we use the machinery of log-structures on toroidal crossing spaces, so we restrict ourselves back to the Gross-Siebert polyhedral base $B$ with elementary simple singularities.

This note consists of two parts. 
The first three sections are devoted to the construction of the compactification of the torus bundle from over the smooth part $B_0$ of the base to all of $B$. The last section compares the topology of the total space of the compactified torus bundle with the Kato-Nakayama space obtained from a toric log Calabi-Yau space.

The primary purpose of this note is an announcement, however, we do give a precise definition of the setup, its basic notions, some discussion of these and the statement of the main results to be achieved. We carry out the compactification construction in dimension three under a unimodularity assumption for illustration. 
Some results may be stated only in special cases and proofs may be sketchy or omitted. All statements in full generality and rigorous proofs will appear soon in \cite{RZ2}.
\begin{acknowledgments}
We are indebted to Mark Gross and Bernd Siebert for sharing their ideas and unpublished notes on the subject, parts of which will enter \cite{RZ2}.
Our gratitude for hospitality goes to Mittag-Leffler Institute, Oberwolfach MFO, University of Miami, MATRIX Institute, JGU Mainz and Kansas State University.
\end{acknowledgments}

\section{Integral affine manifolds with singularities} 
\label{section-mfd-setup}
Let $B$ be a pure $n$-dimensional regular CW complex which is a manifold. We fix the first barycentric subdivision $\bsd B$ of $B$ and let $\bar D$ be the subcomplex of $\bsd B$ which consists of simplices spanned by the barycenters of strata of $B$ which are not vertices and not facets.
That is, $\bar D$ is an $(n-2)$-dimensional subcomplex of $\bsd B$ which lives inside the $(n-1)$-skeleton of $B$ and misses all vertices of $B$.

Suppose that we are given an integral affine structure on $B_{00}:=B\setminus \bar D$. That is, $B_{00}$ is given the structure of a smooth manifold and a flat connection of its tangent bundle $TB_{00}$ with holonomy in $\GL_n(\Z)$. We denote by $\Lambda$ the rank $n$ local system of flat integral vectors in $TB_{00}$. Similar, the local system $\check\Lambda$ stands for the flat integral covectors in the cotangent bundle $T^*B_{00}$.

Each facet of $\bar D$, being of codimension 2, has a small loop around it in $B_{00}$ and we compute the monodromy of the affine structure along this loop. If the monodromy is trivial we can extend the affine structure over this facet. If the monodromy is not trivial, then this facet becomes a part of the true {\bf discriminant} $D$ which is a full-dimensional subcomplex of $\bar D$, that is still a codimension 2 subcomplex of $\bsd B$. We denote by $B_0:=B\setminus D$ the {\bf smooth} part of the base, this is as far as the affine structure extends.

Now we describe the requirements for the monodromy of the affine structure. Let $\iota\colon B_0 \hookrightarrow B$ be the inclusion of the smooth part into the base. Then $\iota_* \Lambda$ and $\iota_* \check\Lambda$ are the constructible sheaves of locally invariant sublattices of $\Lambda$ and $\check\Lambda$. In particular, the stalk of $\iota_* \Lambda$ at a point $x$ in the discriminant $D$ extends as a constant subsheaf of $\Lambda$ in a neighborhood $U$ of $x$ (the $\Lambda$ itself is not trivializable on $U\setminus D$), and similar for $\check\Lambda$.  We denote the restriction of $\iota_* \Lambda$ to $D$ by $L$, and the restriction of $\iota_* \check\Lambda$ to $D$ by $\check L$, both are constructible sheaves on $D$.

Let $x\in D$ be a point which lies in the stratum $\tau$. Pick a nearby base point $y\in B_0$. The local fundamental group of $B_0$ in a neighborhood of $x$ is generated by the loops around the maximal strata of $D$ and we want to see its monodromy image $G_x$ in $\GL( \Lambda_y)$. The minimal requirement is that $G_x$ is an abelian subgroup of $\GL( \Lambda_y)$. In fact we want to require even more.
Since $L,\check L$ are constant on the relative interior of each stratum $\tau$ of $D$, we simply refer to the stalk at any point in that relative interior by $L_\tau$, respectively $\check L_\tau$.

\begin{definition}
For a stratum $\tau\subset D$, suppose there are sublattices $L_1, \dots ,L_r$ in $L_\tau$, linearly independent over $\Q$, and sublattices $\check L_1, \dots, \check L_r$  in $\check L_\tau$, also linearly independent over $\Q$. 
We call the collection of sublattices {\bf semi-simple} if every $L_i$ is orthogonal to every $\check L_j$ (including $i=j$). 
If every stratum $\tau\subset D$ permits a semi-simple collection of sublattices so that the monodromy group $G_x$ for any $x$ in the interior of $\tau$ has the form $\id+L_1\otimes \check L_1 +\dots + L_r\otimes \check L_r$ then we say that $(B,D)$ is an integral affine manifold with {\bf semi-simple abelian} (or for short just semi-simple) singularities.

\end{definition}

We denote the rank of $L_i$ by $k_i$ and the rank of $\check L_i$ by $\check k_i$, and let $\ell:=\sum_i k_i$ and $\check \ell:=\sum_i \check k_i$.
It holds $s:=n-\ell -\check \ell \ge 0$.
The semi-simpleness condition says that in a neighborhood $U$ of $x\in D$ the monodromy matrices in a suitable basis of $\Lambda_y\otimes_\ZZ\QQ$ when acting on column vectors have the shape
$$
\begin{pmatrix}1 & 0 &\cdots & 0 & \boxast & 0 &\cdots & 0\\
0 & 1 &\cdots & 0 & 0 & \ddots &  \ddots & \vdots\\
\vdots & 0 & \ddots & \vdots & \vdots & \ddots &\ddots &0\\
\vdots & \vdots & \ddots & 1 & 0 & \cdots &0 &\boxast\\
0 & \cdots &\cdots & 0 &1 &0 &\cdots&0\\
0 & \cdots &\cdots & 0 & 0 &1 & \ddots&\vdots\\
0 & \cdots &\cdots & 0 &\vdots & \ddots & \ddots & 0\\
0 & \cdots &\cdots & 0 &0 & \cdots & 0 & 1\\
\end{pmatrix}
$$
where the first columns correspond to a basis of $L$, the last rows correspond to a basis of $\check L$ and the $(k_i\times \check k_i)$-size $\boxast$-blocks correspond to the lattices $L_i\otimes \check L_i$.

If $\dim B=3$, then necessarily $r\le 1$ and $D$ is a graph. Following Gross, we call a vertex of $D$ {\bf positive} if $\dim \check L=2$ and we call a vertex {\bf negative} if $\dim L=2$.

In fact we want even more. To every stratum $\tau$ of $D$ we would like to associate two collections of lattice polytopes $(\Delta_1, \dots, \Delta_r)_\tau$ in $L_\tau$ and 
$(\check\Delta_1, \dots, \check\Delta_r)_\tau$ in $\check L_\tau$ such that each $L_i$ is generated by the edge vectors of $\Delta_i$, and similar for $\check L_i$. 
We denote by $\pp$ the collection of $\{\Delta_i, \check\Delta_i\}$ for all strata of $D$. Next we discuss the compatibility of the collection $\pp$ that we require for the inclusion maps $\phi\colon L_\tau \hookrightarrow L_\sigma$ and $\check\phi \colon \check L_\tau \hookrightarrow \check L_\sigma$ for any two incident strata $\tau\prec \sigma$ of $D$. Note that $r_\tau \ge r_\sigma$, and we can always match the number $r$ of polytopes in $\sigma$ and $\tau$  by adding the origins $\{0\}$ to play the role of missing $\Delta, \check\Delta$ to the $\sigma$-collection.

\begin{definition}
The collection of polytopes $\pp$ is {\bf compatible} if for any incident pair $\tau\prec \sigma$ of $D$, after a suitable integral translation, $\phi(\Delta_{i,\tau})$ is a face of $\Delta_{i,\sigma}$ and, similarly after integral translation, $\check\phi(\check\Delta_{i,\tau})$ is a face of $\check\Delta_{i,\sigma}$ for all $i=1,\dots, r$ (up to reordering the indices in $\{1,\dots,r\}$).
\end{definition}

Next, we describe the correlation between $\pp$ and the discriminant $D$. To any polytope $\Delta_i$ one can associate its normal fan. Let $Y_i\subset \R^{k_i}$ be the codimension 1 skeleton of that normal fan. Similarly, $\check Y_i\subset \R^{\check k_i}$ is the codimension 1 skeleton of the normal fan to $\check\Delta_i$. For a point $x\in D$ in a stratum $\tau$ we consider the codimension 2 fans in $\R^n=\R^s\times \R^{k_1}\times\R^{\check k_1}\times\dots \times \R^{k_r}\times\R^{\check k_r} $:
\begin{equation}
\label{eq-def-S}
S_{x,i}:=\R^s \times Y_i\times \check Y_i \times \R^{\ell-k_i+\check \ell-\check k_i}.
\end{equation}
Let $S_x$ be their union: $S_x:= \bigcup_i S_{x,i}$. Note that $\bigcap_i S_{x,i}=\R^s$.
Maximal cones in $S_x$ are labeled by the pairs of edges in $(e,f)$, where $e$ is an edge in $\Delta_i$ and $f$ is an edge in $\check\Delta_i$ for some $i$.


%

\begin{definition}
A compatible collection of polytopes $\pp$ is {\bf normal} if for every point $x\in D$ there is a homeomorphism of its neighborhood $U\subset B$ to an open subset  $V\subset\R^n$ which maps $D\cup U$ to $S_x \cup V$. 
\end{definition}

Finally we make connection between $\pp$ and the monodromy of the affine structure. Let $x$ be a point in $D$ which lies in the stratum $\tau$. Pick a nearby base point $y\in B_0$. We assume that the polytopal collection is compatible and normal. Then the local fundamental group of $B_0$ in a neighborhood of $x$ is generated by the loops around the maximal strata of $D$, which are labeled by the pairs of edges $(e,f)$, $e$ in $\Delta_i$ and $f$ in $\check\Delta_i$, some $i$. Orienting the edges determines an orientation of the loop around the corresponding stratum $\sigma_{e,f}$ of $D$, see [GS] for details.

\begin{definition}
A semi-simple integral affine manifold $(B,D)$ with a compatible normal collection of polytopes $\pp$ is called {\bf semi-simple polytopal} if the local monodromy along the loop $\sigma_{e,f}$ is given by $\id + e \otimes f$. 
\end{definition}

The collection of polytopes $\pp$ is reminiscent of the Batyrev-Borisov nef-partitions. The semi-simple polytopal integral affine manifolds therefore mimic local complete intersections in algebraic geometry. 

We next give an example of a semi-simple affine structure which is not polytopal. The figure below shows a part of discriminant in a 3-dimensional base and the monodromy matrices around the intervals in $D$ with respect to some base point $y\in B\setminus D$ and a suitable basis of $\Lambda_y\cong\ZZ^3$. 
\begin{figure}[h]\label{fig:polytopal}
\begin{tikzpicture}[scale=1]
\draw [thick] (0,0)--(4,0)--(6,1);
\draw [thick] (4,0)--(3.3, -1.3);
\draw [thick] (-.6,1.2)--(0,0)--(-.3,-1.3);
\draw (2,.5) node{$\left(\begin{smallmatrix} 1&0&2\\0&1&0\\0&0&1\end{smallmatrix}\right)$};
\draw (4.5,1) node{$\left(\begin{smallmatrix} 1&0&0\\0&1&-1\\0&0&1\end{smallmatrix}\right)$};
\draw (4.5,-.8) node{$\left(\begin{smallmatrix} 1&0&-2\\0&1&1\\0&0&1\end{smallmatrix}\right)$};
\draw (-1.3,.7) node{$\left(\begin{smallmatrix} 1&-1&0\\0&1&0\\0&0&1\end{smallmatrix}\right)$};
\draw (-1.2,-.8) node{$\left(\begin{smallmatrix} 1&1&-2\\0&1&0\\0&0&1\end{smallmatrix}\right)$};
\draw (2,-.3) node{$\tau$};
\end{tikzpicture}
\caption{A non-polytopal semi-simple affine structure.}
\end{figure}
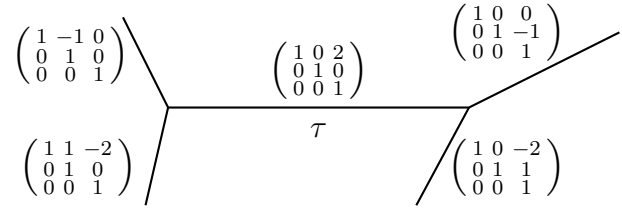
The edge $\tau$ connects the positive vertex on the left with the negative vertex on the right.
The monodromy around the middle edge $\tau$ is twice the standard focus-focus case. 
The point is that it is impossible to decide which of the two intervals $\Delta_\tau$ or $\check\Delta_\tau$ has length 2. 
The vertex on the left requires $\Delta_\tau$ to be length one while the vertex on the right requires $\check\Delta_\tau$ to be length one.

\begin{figure}[h]\label{fig:non-convex}
\begin{tikzpicture}[scale=1]
\draw [thick] (-1.5,1.5)--(0,0)--(2,0);
\draw [thick] (0,-1.6)--(0,0)--(0,1.6);
\draw [thick, ->] (5,-1)--(5,1);
\draw [thick, ->] (5,1)--(6,1);
\draw [thick, ->] (6,1)--(4,-1);
\draw [thick, ->] (4,-1)--(5,-1);
\draw (1.5,-.5) node{$\left(\begin{smallmatrix} 1&0&0\\0&1&2\\0&0&1\end{smallmatrix}\right)$};
\draw (.7,1.1) node{$\left(\begin{smallmatrix} 1&0&1\\0&1&0\\0&0&1\end{smallmatrix}\right)$};
\draw (-2,.9) node{$\left(\begin{smallmatrix} 1&0&-2\\0&1&-2\\0&0&1\end{smallmatrix}\right)$};
\draw (-.7,-1.1) node{$\left(\begin{smallmatrix} 1&0&1\\0&1&0\\0&0&1\end{smallmatrix}\right)$};
\end{tikzpicture}
\caption{Another non-polytopal semi-simple affine structure.}
\end{figure}
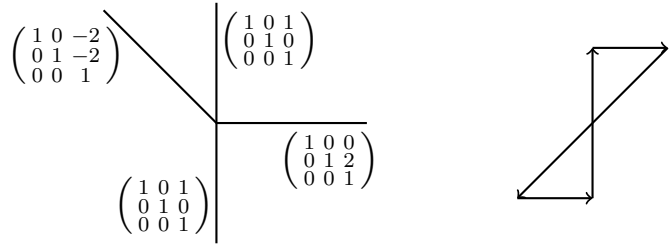

Another feature of a polytopal affine structure is some sort of local convexity of the monodromy. Figure \ref{fig:non-convex} shows an example of a negative vertex in a 3-dimensional base with the monodromy matrices around the four adjacent edges. The monodromy vectors do not form a convex polytope.

%

\section{Local models}
\begin{figure}[htb]
\label{fig:localZ}
$$
\xymatrix{
 &\ar_{\mu_0}[ldd] U_{\Sigma}\ar^(.4){(z^{w_1},\dots, z^{w_r})}[dd]& &\ar[ll] V\ar[d] & \ar_{}[l] \tilde V\ar[d]\\
& & &   (\C^*)^{\ell} \ar_(.5){(f_1,\dots, f_r)}[lld]   \ar^{\text{finite abelian cover}}[d] &\ar_{\Phi}[l] (\C^*)^{\ell} \ar^{\text{finite abelian cover}}[d]\\
\Sigma^\vee \ar^{\mod w_1,\dots, w_r}[d] & \C^r & &\ar^(.65){\bar f_1 \times...\times \bar f_r}[ll](\C^*)^{k_1}\times...\times (\C^*)^{k_r} \ar^{\log}[d] 
     &\ar_{\bar\Phi}[l] (\C^*)^{k_1}\times...\times (\C^*)^{k_r} \ar^{\log}[d]\\
\R^{\check \ell}& & &\R^\ell &\R^\ell
}
$$
\hspace{1.7 in}
\includegraphics[width=3.5in]{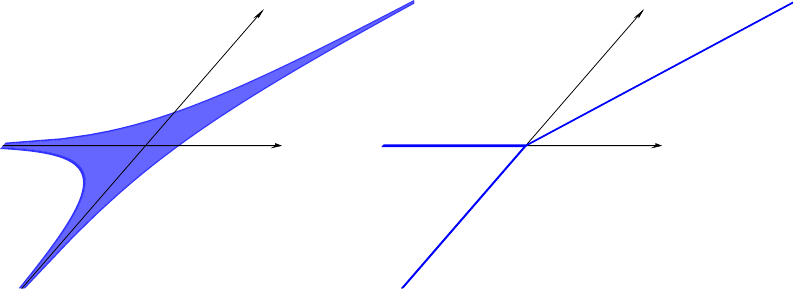}
 \caption{Bottom pictures: $\log (f^{-1}(0))$ and $\log((f\circ \Phi)^{-1}(0))$ for $r=1, \ k=2$.}
 \label{fig:amoeba}
\end{figure}

Let us consider a point $x$ in the interior of a stratum $\tau$ in the discriminant. We will describe a local model of the torus fibration over a neighborhood of $x$ in $B$. 
%
The construction of $X_\tau$ as a fiber product (the left side of the diagram in Figure \ref{fig:localZ}) is pretty standard, see, e.g. [GS]. The novelty here is the rightmost column. 

We explain the details now.
Let $\Sigma$ be the cone over the convex hull $\Conv\{(\check\Delta_i, e_i)\} \subset  \check L_\R \oplus \R^r$, where $e_i=(0,\dots,1,\dots,0)$ is the $i$-th basis vector of $\R^r$, and let $\Sigma^\vee$ be its dual cone and $\Sigma^\vee_\ZZ$ the integral points in the dual cone.
The affine toric variety $U_{\Sigma}=\operatorname{Spec} \C[\Sigma^\vee_\Z]$ has $r$ monomials $z^{w_i}$ corresponding to the integral vectors $w_i$ in $\Sigma^\vee$ defined by 
$$w_i(\check\Delta_i)=1, \quad w_i(\check\Delta_j)=0, \ j\ne i,
$$ 
which gives the map $U_\Sigma \to \C^r$.
The map $\mu_0\colon U_\Sigma \to \Sigma^\vee$ is the moment map, and the cone $\Sigma^\vee$ projects surjectively to $\R^{\check \ell}$ by taking the quotient by the subspace spanned by all $w_i$'s. 

Let $L'\cong \Z^\ell$ be the sublattice in $L$ generated by $L_1,\dots,L_r$ over $\Q$. 
That is $L'=L\cap ((L_1\oplus\dots\oplus L_r)\otimes \Q)$ and $L'$ is the unique direct summand of $L$ that has rank $\ell$ and contains $L_1,...,L_r$. 
The map 
$$ (\C^*)^{\ell} \cong \Hom(L',\CC^*)\ra \Hom(L_1\oplus...\oplus L_r,\CC^*)\cong  (\C^*)^{k_1}\times...\times (\C^*)^{k_r} $$ 
is the abelian cover that appears as the center vertical map in Figure~\ref{fig:localZ} (and a homeomorphic version of it also on the right).
The finite abelian cover takes care of the fact that the sublattice $L_1\oplus\dots\oplus L_r$ may have finite index $>1$ in $L'\subseteq L$. 
There are two ingredients for that. First, the lattice $L_i$ may have an index in its saturation in $L$. 
Second, the direct sum $(L_1 \oplus\dots\oplus L_r)\otimes \Q$ may not split over $\Z$. 

Each polytope $\Delta_i$ defines a function 
$$f_i=\sum_{v\in \operatorname{vert} \Delta_i} c_v z^v \quad :  \quad \Hom(L', \C^*) \to \C, 
\quad \text{for a general choice of} \ c_v\in \C^*.
$$
The same expression also defines a function $\bar f_i\colon \Hom(L_i, \C^*) \to \C$, so that the triangle in Figure~\ref{fig:localZ} commutes. 
The space $V$ is given as the complete intersection $\{f_i=w_i\}$ in $U_\Sigma\times (\C^*)^\ell$.
The full local model for the torus fibration $X_\tau$ over a neighborhood of $x$ is given by multiplying $V$ by the factor of $\log\colon (\C^*)^s \to \R^s$. We reserve the right to split this factor between $U_\Sigma$ and $(\C^*)^\ell$ as needed to match the models for adjacent strata in~$D$.

To actually attach the right column in Figure~\ref{fig:localZ}, we will assume that {\em all polytopes $\Delta_i$ are simplices}. 
At the base of the abelian cover, the map $(f_1,\dots,f_r)\colon (\C^*)^\ell \to \C^r$ splits as a product and, by the assumption of $\Delta_i$ to be a simplex, the hypersurface $\{f_i=0\}$ in $(\C^*)^{k_i}$ is a cover of the pair-of-pants $\{\bar f_i=0\}$. 
The right column is entirely defined once we specify $\bar\Phi=\bar\Phi_1\times\dots\times \bar\Phi_r$ if we additionally require that the two squares adjacent to the right column be Cartesian.
Note that $\Delta_i$ is a unimodular simplex with respect to $L_i$ (by definition of $L_i$).
We will use respectively for $\bar\Phi_i$ the map denoted $\Phi$ in the following key theorem.
For $\Delta$ the standard simplex in $\RR^k$, as before, denote by $Y$ the codimension 1 skeleton of the normal fan to $\Delta$.

\begin{theorem}[\cite{RZ}]\label{thm:ober}
Let $H=\{1+y_1+\dots+y_k=0\} \subset (\C^*)^k$ be the $(k-1)$-dimensional pair-of-pants. Then there is a homeomorphism of the pairs $\Phi\colon ((\C^*)^k,H) \to ((\C^*)^k,\mathcal H)$, where $\mathcal H$ is the ober-tropical pair-of-pants which is mapped by $\log$ to $Y$  with equidimensional fibers.
The homeomorphism restricts well to the boundary
under compactifying $(\C^*)^k$ to the product $\Delta \times \TT^k$ using the moment map $\mu\colon(\C^*)^k \to \Delta^\circ$ (that maps to the interior of the simplex).
\end{theorem} 

The homeomorphism $\Phi$ in the theorem may be viewed as a deformation the $\log$ map so that the image of the pair-of-pants become the tropical hyperplane $Y$ (the spine of the amoeba), rather then the amoeba itself.
Now the fibration $X_\tau \to \R^n$ is induced after applying the homeomophism $\bar\Phi_i$ on the $(\C^*)^{k_i}$ factor (replacing $V$ by a homeomorphic space $\tilde V$). 
The discriminant of the fibration is precisely $S_x= \bigcup_i S_{x,i}$, see \eqref{eq-def-S}, and this has codimension two in $\RR^n$.

Finally, to be able to view the fibration $X_\tau \to \R^n$ as a compactification of the smooth fibration $X_0=T^*B_0/\check\Lambda\ra B_0$ in a neighbourhood of $\tau$, one needs to replace the $\log$ map on the $(\C^*)^\ell$ factor by a suitable (other) moment map $\mu$ so that its image is the interior of a polytope rather than all of $\RR^n$. A straight forward calculation then shows that the monodromy agrees with the local description in the neighborhood of $x\in B$, see, e.g [GS].

\section{Gluing the torus fibration with parameters}
Let $B$ be an integral affine manifold with semi-simple polytopal singularities so that each $\Delta_i$ is a simplex.
\begin{theorem}[\cite{RZ2}]
\label{main-thm1}
There is a topological orbifold $X$ which compactifies the torus bundle $X_0=T^*B_0/\check\Lambda$ to a fibration $X\to B$ with $n$-dimensional fibers (singular over $D\subset B$). If all local cones $\Sigma_\tau$ are unimodular simplicial cones then $X$ is a manifold.
\end{theorem}

In fact one can vary the gluing data (a.k.a.~B-field) to get a whole family of torus fibrations $X_\gamma$ over~$B$. 
The parameter space of gluings is a torsor over $H^1(B,\iota_*\check\Lambda \otimes U(1))$. 
One can make sense of the parameter space itself being $H^1(B,\iota_*\check\Lambda \otimes U(1))$ by carefully choosing preferred matching sections for local models - for these additional constructions, we refer to \cite{RZ2}. 
By the universal coefficient theorem, $H^1(B,\iota_*\check\Lambda)\otimes_\ZZ U(1)\subseteq H^1(B,\iota_*\check\Lambda \otimes U(1))$ is the component of the identity.
We will relate the resulting family
\begin{equation}
\label{eq-top-family}
\shX\ra H^1(B,\iota_*\check\Lambda)\otimes_\ZZ U(1)
\end{equation}
to the Gross-Siebert program in the next section. 

For the remainder of this section, we are going to carry out the compactification procedure in dimension 3 (and thus for $r=1$ at all strata of $D$ as we already pointed out in Section~\ref{section-mfd-setup}) and we additionally assume that both $\Delta_\tau$ and $\check\Delta_\tau$ are unimodular simplices. We follow the approach of Gross \cite{Gross}, that is, we successively compactify the fibration over the star neighborhoods of vertices in $D$ (the barycenters of faces in $B$) ordered by dimension of the corresponding face in $B$.

We begin with $X_0=T^*B_0/\check\Lambda$.
As mentioned before, there are two types of vertices in $D$: barycenters of 2-faces (negative vertices if more than bivalent) and barycenters of 1-faces (positive vertices if more than bivalent) in $B$. We will first compactify over the star neighborhoods (in the $\bsd B$) of the 2-face vertices.
Let $x\in D$ be the barycenter of a 2-face $Q$ of $B$. Then $\check\Delta$ is the unit interval and there are two possibilities for $\Delta_x$: the standard 2-simplex or the unit interval.

Case 1: For $\Delta$ the unit interval, we have the $(2,2)$-case in \cite{Gross} which is the standard focus-focus compactification times a $\C^*$-factor.

Case 2: For $\Delta$ the standard 2-simplex, we find $D$ has a trivalent vertex at $x$ and this is referred to as the $(2,1)$-case in \cite{Gross}.  Let $\Delta^\circ$ be the interior of $\Delta$ and $Y$ be the union of the 3 intervals in $\bsd \Delta$ that connect the barycenter of~$\Delta$ with the barycenter of an edge of $\Delta$ respectively.
In a neighborhood of $x$, the torus bundle $T^*B_0/\check\Lambda$ becomes a trivial $\TT^2$-bundle once we take the quotient by the coinvariant (vanishing) circle $(\iota_* \check\Lambda)_x \otimes U(1)$, hence this $\TT^2$-bundle extends over the discriminant. Precisely, to glue in the cotangent torus bundle $T^*\Delta^\circ/ L^* \cong \Delta^\circ \times \TT^2 \cong Q\times (\check\Lambda/\check L\otimes U(1)) $ over the simplex (here $L^*$ is the dual lattice to $L$) we just need to identify the interior of the simplex $\Delta^\circ$ with the cell~$Q$. 
Let $Q_1, Q_2, Q_3$ be the 3 boundary intervals of $Q$ which meet $D$.

\begin{lemma}\label{lemma:S^1}
There is a homeomorphism $\psi\colon (\Delta^\circ, Y) \to (Q, D\cap Q)$ which extends to a homeomorphism including the 3 boundary intervals of $\Delta$  and identifying these with the 3 boundary intervals $Q_1, Q_2, Q_3$ of~$Q$.
\end{lemma}

Thus we have a well-defined (trivial)  $\TT^2$-bundle (the quotient by the $\check L$-circle) over the star neighborhood of $x$ in $\bsd B$.
Our next step is to compactify the circle bundle over the 5-dimensional manifold $((\Delta^\circ \times \TT^2) \times \R)  \setminus  (\mathcal H \times\{0\})$ to a fibration over $\Delta^\circ \times \TT^2\times\R$. 
Here, $\mathcal H\subset (\C^*)^2 \cong \Delta^\circ \times \TT^2$ denotes the ober-tropical pair-of-pants from \cite{RZ}, a 2-dimensional submanifold of $\Delta^\circ \times \TT^2\times\R$. 
We can either glue in the local model $X_Q$ from the previous section, or use the following proposition, leading to a homeomorphic result:

\begin{lemma}[cf. \cite{Gross}, Proposition 2.5] \label{lemma:21}
Let $U$ be the complement of an oriented connected submanifold $S$ of codimension 3 in a manifold $\bar U$ and let $\pi\colon X\to U$ be a principal $S^1$-bundle with the Chern class $c_1=\pm\kappa$ in $H^3_S(\bar U, \Z) \cong H^0(S,\Z)$ for some $\kappa>0$. Then there is a unique compactification  to an orbifold $\bar X=X \cup S$ such that $\bar \pi\colon \bar X \to \bar U$ is a proper map and $\bar X$ is a manifold if $\kappa=1$.
\end{lemma}
Our local description of the monodromy implies that $\kappa$ coincides with the lattice length of $\check\Delta_x$ which we assumed to be one for this article. As pointed out by \cite{Gross} already, changing the orientation of the $S^1$-action, changes the sign of $c_1$.
The label ``negative'' for the vertex $x$ of $D$ doesn't stem from the sign of $c_1$ but the Euler number of the local model $X_Q$ which is $-1$ (and the same as the Euler number of the fiber over $x$).

\begin{figure}[h]\label{fig:wiggle}
\begin{tikzpicture}[scale=1]

\draw  (0,0)--(6,0)--(0,6)--(0,0);
\draw [thick, blue] (2,2)--(3,0);
\draw [thick, blue] (2,2)--(3,3);
\draw [thick, blue] (2,2)--(0,3);
\draw [dashed, fill=blue!30, opacity=.5] (1,4)--(0.66,3.66)--(1,2.5)--(.66, 1.66) --(1,1)--(1.66,0.66)--(2.5, 1)--(3.66,0.66)--(4,1)--(3.66,1.66)--(2.5, 2.5)--(1.66,3.66)--cycle;

\draw (-1.5,3.5)--(-1,3.5)--(-1,4)--(-1.5,4)--(-1.5,3.5);
\draw [thick, red] (-1.5, 3.75)--(-1.33, 3.83)--(-1.16, 3.66)--(-1,3.75);
\draw [thick, red] (-1.34, 3.83)--(-1.25,4);
\draw [thick, red] (-1.16, 3.66)--(-1.25,3.5);
\draw [->] (-1,3.75) to [out=0 ,in=100] (2,2.02);

\draw (3.5,4.5)--(4,4.5)--(4,5)--(3.5,5)--(3.5,4.5);
\draw [thick, red] (3.5, 4.75)--(3.66, 4.83);
\draw [thick, red] (3.83, 4.66)--(4,4.75);
\draw [thick, red] (3.66, 4.83)--(3.75,5);
\draw [thick, red] (3.83, 4.66)--(3.75,4.5);
\draw [->] (3.75,4.5) to [out=220 ,in=120] (2.25,2.25);

\draw (4.5,3.5)--(5,3.5)--(5,4)--(4.5,4)--(4.5,3.5);
\draw [thick, red] (4.5, 3.75)--(4.75,4);
\draw [thick, red] (4.75,3.5)--(5,3.75);
\draw [->] (4.75,3.5) to [out=220 ,in=-30] (2.75,2.75);

\draw (5.5,2.5)--(6,2.5)--(6,3)--(5.5,3)--(5.5,2.5);
\draw [thick, red] (5.75,2.5)--(5.83, 2.66)--(5.66, 2.83)--(5.75,3);
\draw [->] (5.5,2.75) to [out=180 ,in=0] (2.25,1.5);

\draw (6.5,1.5)--(7,1.5)--(7,2)--(6.5,2)--(6.5,1.5);
\draw [thick, red] (6.75,1.5)--(6.75,2);
\draw [->] (6.5,1.75) to [out=180 ,in=0] (2.75,.5);

\draw (-1.5,1.5)--(-1,1.5)--(-1,2)--(-1.5,2)--(-1.5,1.5);
\draw [thick, red] (-1.5, 2.75)--(-1,2.75);
\draw [->] (-1,1.75) to [out=0 ,in=240] (1.5,2.25);

\draw (-1.5,2.5)--(-1,2.5)--(-1,3)--(-1.5,3)--(-1.5,2.5);
\draw [thick, red] (-1.5, 1.75)--(-1.33, 1.83)--(-1.16, 1.66)--(-1,1.75);
\draw [->] (-1,2.75) to [out=-30 ,in=240] (.5,2.75);

\end{tikzpicture}
\caption{Fading off the wiggling of red circles along $Y\subset \Delta^\circ$.}
\end{figure}
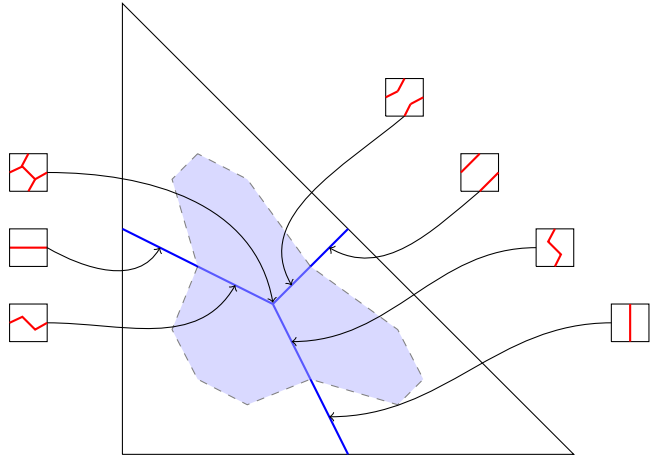

Now comes the important step of {\bf unwiggling} the ober-tropical fiber tori for extending the compactification over the rest of $D$. Figure \ref{fig:wiggle} shows the $\TT^2$-fibers over different points in $\Delta^\circ$ indicated by little squares. 
The red locus inside each square is the intersection of $\mathcal H$ with the corresponding $\TT^2$-fiber, the ``ober-tropical fibers'' over points of $Y$.
As we move away from $x\in D$ we deform the ober-tropical fibers so that they become more and more straight circles. Outside the second barycentric star of $x$ (the shaded region) they are true linear circles in~$\TT^2$ and are ready to be glued with the neighboring model. The $S^1$-fibration over $\Delta^\circ \times \TT^2\times\R$ collapses precisely over the red circles in the $\T^2$-fibers over $Y\times \{0\} \subset \Delta^\circ\times \R$.

The main point of unwiggling is to achieve the following property: close to the boundary of the cell $Q$ the constructed space $X$ may not only be thought not only as an $S^1$-fibration over $Q \times\R \times \TT^2$ but also as a $\TT^2$-fibration over $Q \times\R \times S^1$ via taking the quotient by the circle in the base $\TT^2$ that is the homotopy class of the respective red ober-tropical fiber circle.  From this perspective, the fibers over $D\times S^1$ are pinched tori (homeomorphic to $I_1$-degenerate elliptic curves). This helps us to do the last step, namely compactify the fibration over the vertices of $D$ which are barycenters of the 1-dimensional strata in $B$.

Let us finally discuss the compactification over a vertex $x\in D$ that is the barycenter of a one-cell in $B$. 
At this barycenter, $\Delta$ is the unit interval. If $\check\Delta$ is also a unit interval then we are back to the focus-focus $(2,2)$-case which is straightforward to compactify, so assume $\check\Delta$ is a standard 2-simplex, this is the $(1,2)$-case in \cite{Gross}.

We state a more general result that relates back to Figure~\ref{fig:localZ}. Let $\Sigma$ be a cone in $\R^k \times \R$ over a lattice simplex $\check\Delta \subset \R^k \times \{1\}$, and let $\Sigma^\vee$ be the dual cone. The projection $\R^k \times \R \to \R$ gives a linear map of cones $w\colon \Sigma \to \R_{\ge 0}$ which, in turn, defines a map of affine toric varieties $z^w\colon U_\Sigma = \operatorname{Spec}\C[\Sigma^\vee_\Z] \to \C$.  Let $\mu_0\colon U_\Sigma \to \R^k$ be the moment map with respect to the $\TT^k$-action on $U_\Sigma$ which fixes $z^w$.

We consider the $\TT^k$-torus fibration $\pi\colon U_\Sigma \to \R^k\times \C$ given by $z\mapsto (\mu_0, z^w)$. Recall that $\check Y\subset \R^k$ stands for the codimension 1 skeleton of the normal fan to $\check \Delta$. Away from $\check Y \times \{0\}$, the map $\pi$ is a $\TT^k$-bundle with the monodromy prescribed by the pair $(\Delta=[0,1],\check\Delta)$. 
Over the strata of $\check Y$, the fibers of $\pi$ become lower-dimensional tori (reflecting the dimension of the stratum) with the fiber over $\{0\}$ being just a single point.
We denote by $\pi_0: U_\Sigma\setminus \{0\} \to (\R^k\times\C) \setminus \{(0,0)\}$ the restriction of the fibration $\pi$ to the complement of the origin.

\begin{lemma}[cf. \cite{Gross}, Proposition 2.9, for the 3-dimensional case] \label{lemma:12}
 Let $X \to (\R^k\times\C) \setminus \{(0,0)\}$ be a torus fibration homeomorphic to $\pi_0$.
There is a unique one point compactification to an orbifold $\bar X=X \cup \{pt\}$ such that $\bar \pi_0\colon \bar X \to  \R^k\times\C$ is a proper map. 
Consequently, the fibration $\bar \pi_0\colon \bar X \to  \R^k\times\C$ is homeomorphic to $\pi\colon U_\Sigma \to \R^k\times \C$ and
$\bar X$ is a manifold if $\check\Delta$ is a unimodular simplex.
\end{lemma}

There is a generalization of this statement when the pair $(\C,\{0\})$ is replaced by a pair $(U,S)$ of $S$ being a submanifold in $U$ of codimension 2. This, in particular, covers Lemma~\ref{lemma:S^1} as a special case $k=1$. The relevant case for us is $k=2$, $U=\R\times (\RR/\ZZ)$ and $S$ is the point $(0,0)\in \R\times (\RR/\ZZ)$. We identify $S^1=\RR/\ZZ$ and refer to $0$ as the corresponding point in $S^1$ in the following.

First, we note that similar to the $(2,1)$-vertex, the torus bundle $T^*B_0/\check\Lambda$ in a neighborhood $W$ of $x$ becomes a trivial $S^1$-bundle once we take the quotient by the coinvariant (vanishing) $\TT^2$-subbundle $(\iota_* \check\Lambda)_x \otimes S^1$, thus it extends over $D$. Second, we may view  the torus bundle $T^*(B_0\cap W)/\check\Lambda$ as a $\TT^2$-bundle over $(W\setminus D)\times S^1$.

\begin{lemma}
The $\TT^2$-bundle over $(W\setminus D)\times S^1$ extends to a singular $\TT^2$-fibration $\pi_0\colon X\to (W\times S^1) \setminus (x \times \{0\})$ by adding the $I_1$-fibers over $(\check Y \setminus x) \times S^1$. The resulting fibration $X$ agrees with those coming from the neighboring $(2,1)$-vertices of $D$ after the unwiggling. 
Moreover, the fibration $X\to (W\times S^1) \setminus (x \times \{0\})$ is homeomorphic to $\pi_0$, so satisfies the hypothesis of Lemma~\ref{lemma:12}.
\end{lemma}
Applying Lemma~\ref{lemma:12}, the space $X$ compactifies to $\bar X \to W\times S^1$ by adding the point $x \times \{0\}$. This completes the compactification process.

Finally we briefly comment on the parameter space $H^1(B,\iota_*\check\Lambda \otimes U(1))$ of gluings. Already in building the smooth part $T^*B_0/\check\Lambda$ one can twist by a \v{C}ech cocycle representing an element in $H^1(B_0, \check\Lambda \otimes U(1))$. Furthermore when gluing in the local models around the vertices the twisting can be made when identifying the $\TT^3$-fibers of $T^*B_0/\check\Lambda$ with the $S^1$-bundle over $\TT^2$ (for the $(2,1)$-vertices) or $\TT^2$-bundles over $S^1$ (for the $(1,2)$-vertices). Lastly, when identifying the models between $(2,1)$ and $(1,2)$-vertices there is only $\TT^2$-freedom of twistings which corresponds to the sheaf $\iota_*\check\Lambda$ dropping the rank along the edges of $D$.

\section{Canonical Calabi-Yau families and their Kato-Nakayama spaces}
Recall from \cite[Theorem 5.2, Theorem 5.4]{logmirror1} and \cite[Remark 5.3]{logmirror2} that, given an integral affine manifold with simple singularities $B$ and a compatible polyhedral decomposition $\P$ with multivalued strictly convex piecewise affine function $\varphi$ with integral slopes, there is an associated algebraic family of toric log Calabi-Yau spaces\footnote{We have implicitly picked a splitting of the surjection $H^1(B,\iota_*\check\Lambda)\ra H^1(B,\iota_*\check\Lambda)/H^1(B,\iota_*\check\Lambda)_{\tors}$.}
\begin{equation} \label{GS-family}
X_0(B,\P,\varphi)\ra  S:=\Spec \CC[H^1(B,\iota_*\check\Lambda)^*]
\end{equation} 
which is semi-universal by \cite[Theorem C.6]{RS}.
Analytification and application of the Kato-Nakayama functor associates to the log morphism \eqref{GS-family} a continuous surjection of topological spaces 
$$\shX''\ra \shS'':=\big( H^1(B,\iota_*\check\Lambda)\otimes_\ZZ\CC^*\big) \times U(1)$$
which is a fiber bundle by \cite[Theorem 5.1]{NO}.
We restrict the family to 
\begin{equation} \label{U1-family}
\shX'\ra \shS':=\big( H^1(B,\iota_*\check\Lambda)\otimes_\ZZ U(1)\big) \times U(1).
\end{equation}
Let $c_1(\varphi)$ denote the class of $\varphi$ in $H^1(B,\iota_*\check\Lambda)$.
The inclusion $\ZZ c_1(\varphi)\subseteq  H^1(B,\iota_*\check\Lambda)$ induces a map of real Lie groups
$$\phi_{\varphi}\colon U(1)\ra H^1(B,\iota_*\check\Lambda)\otimes_\ZZ U(1) $$
where we have identified $(\ZZ c_1(\varphi))\otimes U(1)=U(1)$.

As explained in \cite[\S4.1]{RS}, there is an equivariant $U(1)$ action on the family \eqref{U1-family} and for the base space, by \cite[(4.14)]{RS}, it is given by
$$U(1)\times \shS'\ra \shS',\quad \lambda.(s,t)= \big(\phi_{\varphi}(\lambda)\cdot s,\lambda^{-1}\cdot t\big).$$
Consequently, the family \eqref{U1-family} is a base change of the restricted family $\shX\ra \shS:= H^1(B,\iota_*\check\Lambda)\otimes_\ZZ U(1)$ under the base change homomorphism 
$$ \id\times\phi_{\varphi}\colon  \big( H^1(B,\iota_*\check\Lambda)\otimes_\ZZ U(1)\big) \times U(1)\ra H^1(B,\iota_*\check\Lambda)\otimes_\ZZ U(1).$$
The relevant topological information is therefore already contained in $\shX\ra \shS$. 
In \cite[Section 2.1]{RS}, a moment map $X_0(B,\P,\varphi)\ra B$ was given under the assumption that $X_0(B,\P,\varphi)$ is projective over $ S$. Since we restricted to the $U(1)$-part of the gluing torus when taking \eqref{U1-family}, a moment map exists even without the projectivity condition. Composing with the log forget morphism yields a fibration
$$\pi: \shX\ra B$$
whose discriminant has codimension one (being a union of amoebae in real hyperplanes). In our upcoming work \cite{RZ2}, we are going to prove the following result which may be viewed as deforming $\pi$ near the discriminant so that the new discriminant has codimension two.
In other words, up to this deformation, the Gross-Siebert fibration $\pi$ agrees with the compactification of the Strominger-Yau-Zaslow fibration constructed in the previous sections.

\begin{theorem}[\cite{RZ2}] 
The topological family $\shX\ra \shS$ is homeomorphic to the family of compactified torus bundles given in \eqref{eq-top-family} (over the identity on $\shS$).
\end{theorem}
The homomorphism in the theorem commutes with the respective torus fibration maps to $B$ away from a tubular neighbourhood of the discriminant $D\subset B$.

\end{document}